\documentclass{amsart}

\newtheorem{thm}{Theorem}

\newtheorem{lem}[thm]{Lemma}
\newtheorem{cor}[thm]{Corollary}

\newtheorem{prop}[thm]{Proposition}

\newtheorem{conj}[thm]{Conjecture}
   
\theoremstyle{definition}
\newtheorem{defn}[thm]{Definition}

\newtheorem{say}[thm]{}
\newtheorem{exmp}[thm]{Example}

\newtheorem{rem}[thm]{Remark}          
\newtheorem{note}[thm]{Note}            
\newtheorem{ack}{Acknowledgments}        
\newtheorem{notation}[thm]{Notation}   
  
\newtheorem{defn-thm}[thm]{Definition--Theorem}  

\theoremstyle{remark}
\newtheorem{claim}[thm]{Claim}

\newtheorem{outline}[thm]{Outline of the proof}

\renewcommand{\c}[0]{{\mathbb C}}  

\renewcommand{\o}[0]{{\mathcal O}} 
\newcommand{\z}[0]{{\mathbb Z}}

\newcommand{\p}[0]{{\mathbb P}}

\newcommand{\q}[0]{{\mathbb Q}}

\newcommand{\qtq}[1]{\quad\mbox{#1}\quad}

\newcommand{\mult}[0]{\operatorname{mult}}

\newcommand{\aut}[0]{\operatorname{Aut}}

\newcommand{\ric}[0]{\operatorname{Ricci}}

\newcommand{\lcm}[0]{\operatorname{lcm}}
\newcommand{\wdeg}[0]{\operatorname{wdeg}}

\def\fract#1#2{\raise4pt\hbox{$ #1 \atop #2 $}}
\def\decdnar#1{\phantom{\hbox{$\scriptstyle{#1}$}}
\left\downarrow\vbox{\vskip15pt\hbox{$\scriptstyle{#1}$}}\right.}

\def\bfa{{\bf a}}
\def\bfb{{\bf b}}

\def\bfw{{\bf w}}
\def\bfx{{\bf x}}

\def\bfz{{\bf z}}

\def\cals{{\mathcal S}}

\def\bbc{{\mathbb C}}

\def\bbp{{\mathbb P}}
\def\bbq{{\mathbb Q}}

\def\bbz{{\mathbb Z}}

\def\gro{\omega}

\def\grt{\tau}

\def\gsp1{{\mathfrak s}{\mathfrak p}(1)}

\def\ker{\hbox{ker}}

\def\la#1{\hbox to #1pc{\leftarrowfill}}
\def\ra#1{\hbox to #1pc{\rightarrowfill}}
\def\Se{Sasakian-Einstein }
\def\Ke{K\"ahler-Einstein }

\begin{document}
\bibliographystyle{amsalpha}

\title{Einstein Metrics on Spheres}
\author{Charles P. Boyer, Krzysztof Galicki and  J\'anos Koll\'ar}
\address{CPB and KG: Department of Mathematics and Statistics, 
University of New Mexico,
Albuquerque, NM 87131.}
\email{cboyer@math.unm.edu}
\email{galicki@math.unm.edu}
\address{JK: Department of Mathematics,
Princeton University,
Princeton, NJ 08544-1000.}
\email{kollar@math.princeton.edu}

\maketitle

\section{Introduction} 

Any sphere $S^n$ admits a metric of constant sectional
curvature. These canonical metrics are homogeneous and Einstein, that is
the Ricci curvature is a constant multiple of the metric.
The spheres $S^{4m+3}$, $m>1$ are known to have
another $Sp(m+1)$-homogeneous Einstein metric discovered by
Jensen \cite{MR50:5694}. 
In addition, $S^{15}$ has a third ${\rm  Spin}(9)$-invariant
homogeneous Einstein metric discovered by Bourguignon and Karcher
\cite{MR58:12829}. In 1982 Ziller proved that these are the only
homogeneous Einstein metrics on spheres \cite{MR84h:53062}. 
No other Einstein metrics on spheres were known until
1998 when B\"ohm constructed infinite sequences of 
non-isometric Einstein metrics, of positive scalar curvature, on 
$S^5$, $S^6$, $S^7$, $S^8$, and $S^9$ \cite{MR99i:53046}. B\"ohm's metrics are
of cohomogeneity one and they are not only the first 
inhomogeneous Einstein metrics on spheres but also the first
non-canonical Einstein metrics on even-dimensional spheres.
Even with B\"ohm's result, Einstein metrics on spheres appeared to
be rare.

The aim of this paper is to demonstrate that on the contrary, at least
on odd-dimensional spheres, such metrics occur with 
abundance in every dimension. Just as in the case of B\"ohm's construction, 
ours are only existence results. However, we also answer   in the affirmative 
the long standing open question about the existence of Einstein metrics
 on exotic spheres. 
These are differentiable manifolds that are homeomorphic but not
diffeomorphic to a standard sphere $S^n$.

Our method proceeds  as follows. For a sequence
$\bfa=(a_1,\ldots,a_m)\in\z_+^m$  consider the 
{\it Brieskorn--Pham singularity}
$$
Y(\bfa):=\bigl\{\sum_{i=1}^mz_i^{a_i}=0\bigr\}\subset \c^m
\qtq{and its link}
L(\bfa):=Y(\bfa)\cap S^{2m-1}(1).
$$
$L(\bfa)$ is a smooth, compact,  $(2m-3)$-dimensional manifold.
$Y(\bfa)$ has a natural $\c^*$-action and $L(\bfa)$
a natural $S^1$-action (cf. \S \ref{BP.defn}).
When the sequence $\bfa$  satisfies 
certain numerical conditions, we use the  continuity method to produce
an orbifold  K\"ahler-Einstein metric on the quotient $(Y(\bfa)\setminus \{0\})/\c^*$ 
which then
can be lifted to an Einstein metric on the link $L(\bfa)$.
We get in fact more:
\begin{enumerate}
\item[$\bullet$] The connected component of the isometry group of the  metric
is $S^1$.
\item[$\bullet$]
  We  construct continuous  families of
inequivalent Einstein metrics.
\item[$\bullet$]
The K\"ahler-Einstein structure on the quotient $(Y(\bfa)\setminus \{0\})/\c^*$
lifts to a
{\it Sasakian-Einstein metric} on $L(\bfa)$. Hence, these
metrics have  {\it real Killing spinors}
 \cite{FK90} which play an important role 
in the context of p-brane solutions 
in superstring theory and  in M-theory. See also
\cite{ghp} for related work.
\end{enumerate}

In each fixed dimension $(2m-3)$  we obtain a K\"ahler-Einstein metric
on  infinitely many different  quotients $(Y(\bfa)\setminus \{0\})/\c^*$, but
 the link $L(\bfa)$ is  a homotopy sphere 
only for 
finitely many of them. 
 Both the number of inequivalent families of \Se metrics 
and the dimension of their moduli grow double exponentially 
with the dimension.

There is nothing
special about restricting to spheres or even to Brieskorn-Pham type --
our construction is far more general. All the restrictions
made in this article are very far from being optimal  
and we hope that many more cases will be settled in the future.
Even with the current weak conditions 
we get an abundance of new Einstein metrics. 

\begin{thm} On $S^5$ we obtain 68 inequivalent families of 
\Se metrics. Some
of these admit  non-trivial continuous \Se deformations.
\end{thm}

The biggest family, constructed in
Example \ref{12.dim.on 5} has (real) dimension 10.

The metrics we construct  are almost always inequivalent
not just as Sasakian structures
 but also as Riemannian metrics. The
only exception is that a hypersurface  and its conjugate
lead to isometric Riemannian metrics,  see \S\ref{isom.sasak}.

In the next odd dimension the situation becomes
much more interesting. An easy computer search  finds 
8,610 distinct families of \Se structures on standard and exotic 7-spheres.
By Kervaire and Milnor there are 28 oriented
diffeomorphism types of topological 7-spheres \cite{ker-mil}.
(15 types if we ignore orientation.) The results of
Brieskorn allow one to decide which $L(\bfa)$ corresponds to
which exotic sphere \cite{briesk}.
 We get:

\begin{thm}  All 28 oriented diffeomorphism classes on $S^7$ admit
inequivalent families of \Se structures.
\end{thm}

In each case, the number of families is  easily computed and they range
from $231$ to $452$, see \cite{bgkt} for the computations.
 Moreover, there are fairly large moduli. 
For example, the standard 
7-sphere admits an 82-dimensional family of \Se metrics,
 see Example~\ref{12.dim.on 5}. Let us mention here that any 
orientation reversing 
diffeomorphism takes a \Se metric into an Einstein 
metric, but not necessarily a \Se metric, since the 
Sasakian structure fixes the orientation.

Since Milnor's discovery of exotic spheres \cite{MR18:498d}
the study of special 
Riemannian metrics on them has always attracted a lot of attention.
Perhaps the most intriguing question is whether exotic spheres
admit metrics of positive sectional curvature. This problem remains
open. In 1974 Gromoll
and Meyer wrote down a metric of non-negative sectional curvature
on one of the Milnor  spheres \cite{MR51:11347}. 
More recently it has been observed
by Grove and Ziller that all exotic $7$-spheres which are
$S^3$ bundles over $S^4$ admit  metrics of
non-negative sectional curvature \cite{MR2001i:53047}. But it is not known if
any of these metrics can be deformed to a metric of strictly positive 
curvature. Another interesting question concerns the existence
of metrics of positive Ricci curvature on exotic 7-spheres.
This question has now been settled by the result of Wraith who
proved that all spheres that are boundaries of parallelizable 
manifolds admit a metric of positive Ricci curvature \cite{MR98i:53058}. 
A proof of this result 
using techniques similar to the present paper was recently given in
 \cite{BGN03b}.
In dimension $7$
all homotopy spheres have this property. In this context the result of
Theorem 2 can be rephrased as to say that all homotopy $7$-spheres admit
metrics with positive constant Ricci curvature. Lastly,
we should add that although heretofore
 it was unknown whether Einstein metrics existed on 
exotic spheres, Wang--Ziller, Kotschick and Braungardt--Kotschick
studied  Einstein
metrics  on manifolds which are homeomorphic
but not diffeomorphic \cite{ wa-zi, kot, b-kot}.
In dimension 7 there are even examples of homogeneous
Einstein metrics with this property \cite{MR89c:57042}. Kreck and Stolz
find that there are $7$-dimensional manifolds with the 
maximal number of $28$ smooth structures, each of which admits an
Einstein metric with positive scalar curvature. Our Theorem 2 establishes
the same result for $7$-spheres.

In order to organize the higher dimensional cases, note that
every link $L(\bfa)$ bounds a parallelizable manifold
(called the Milnor fiber).
Homotopy $n$-spheres that bound a parallelizable manifold
form a group, called the {\it Kervaire-Milnor group}, denoted by  $bP_{n+1}$.
When $n\equiv 1\mod 4$
the Kervaire-Milnor group  has at most 2 elements,
the standard sphere and the {\it Kervaire sphere}. 
(It is not completely understood in which dimensions are they different.)

\begin{thm} For $n\geq 2$,  the $(4n+1)$-dimensional 
standard and Kervaire spheres 
 both admit many families of inequivalent 
\Se metrics.
\end{thm}

A partial computer search yielded more than $3\cdot 10^6$
cases for $S^9$ and more than $10^9$ cases for $S^{13}$,
including a 21300113901610-dimensional  family, see 
Example~\ref{big.moduli.exmp}.
The only Einstein metric on $S^{13}$ known previously
was the standard one.

In the remaining case of $n\equiv 3\mod 4$ the situation is
 more complicated.  For these values of $n$
the group $bP_{n+1}$ is quite large (see \S\ref{bP2n.defn})
 and we do not know how to show that
every member of it admits a \Se structure,
since our methods do not apply to the examples given in
\cite{briesk}.  We believe, however, that this is
true:

\begin{conj} All odd-dimensional homotopy spheres which
bound parallelizable manifolds admit 
\Se metrics.
\end{conj}

This was checked by computer  in dimensions
up to  15 \cite{bgkt}.

\begin{outline}
Our construction can be divided into four main steps, each of
quite different character.
The first step, dating back to Kobayashi's circle bundle construction
 \cite{Kob63}, is to observe 
that a positive K\"ahler-Einstein metric on the base space of a circle
 bundle gives an Einstein 
metric on the total space. This result was generalized to orbifolds giving 
\Se metrics in 
\cite{BG00}. Thus, a positive K\"ahler-Einstein orbifold metric on 
$(Y(\bfa)\setminus\{0\})/\c^*$
yields  a \Se metric on 
$L(\bfa)$.
In contrast with the cases studied in \cite{BG01, BGN03},
our quotients are not well formed, that is, some group 
elements have codimension 1 fixed point sets.

The second step is to use the continuity method developed by
\cite{aub, siu1, siu2, tian} to construct  K\"ahler-Einstein
metrics on orbifolds. With minor modifications,
the method of \cite{nad, dem-koll} arrives at a sufficient
 condition, involving the
integrability of inverses  of polynomials
on $Y(\bfa)$. These kind of orbifold metrics were first used
in \cite{ti-ya}.

The third step is to check these  conditions.
Reworking the earlier estimates given in  \cite{jk1, BGN03}
already gives some examples, but here we also give an improvement.
This is still, however, quite far from what one would expect.

The final step is to get examples, partly through computer searches,
partly through writing down well chosen sequences.
The closely related  exceptional singularities of \cite{is-pr}
all satisfy our conditions.
\end{outline}

\section{Orbifolds as quotients by $\c^*$-actions}

\begin{defn}[Orbifolds]\label{orbif.defn}
  An {\it orbifold} is 
 a normal, compact,  complex space $X$
 locally given by charts
written as  quotients of smooth coordinate charts.
That is, $X$ can be covered by
open charts $X=\cup U_i$ and for each $U_i$ there
is a smooth complex space $V_i$ and a finite group
$G_i$ acting on $V_i$ such that $U_i$ is
biholomorphic to the quotient space $V_i/G_i$.
The quotient maps are denoted by $\phi_i:V_i\to U_i$.

The {\it classical} (or well formed) case is when 
the fixed point set of every non-identity element of every
$G_i$ has codimension at least 2.
In this case $X$ alone determines the orbifold structure.

One has to be more careful when there are 
codimension 1 fixed point sets. (This happens to be the case 
in all our examples leading to Einstein metrics.)
 Then the quotient map $\phi_i:V_i\to U_i$ has 
branch divisors $D_{ij}\subset U_i$
and ramification divisors $R_{ij}\subset V_i$.
Let $m_{ij}$ denote the ramification index
over $D_{ij}$. Locally near a general point of
$R_{ij}$ the map $\phi_i$ looks like
$$
\c^n\to \c^n,\quad \phi_i:(x_1,x_2,\dots,x_n)\mapsto 
(z_1=x_1^{m_{ij}},z_2=x_2,\dots,z_n=x_n).
$$
Note that
$$
\phi_i^*(dz_1\wedge \dots\wedge dz_n)=
m_{ij}x_1^{m_{ij}-1}\cdot dx_1\wedge \dots\wedge dx_n.
\eqno{\ref{orbif.defn}.1}
$$
The compatibility condition between the charts
that one needs to assume is that  there are
global divisors $D_j\subset X$ and ramification
indices $m_j$ such that $D_{ij}=U_i\cap D_j$
and $m_{ij}=m_j$ (after suitable re-indexing).

It will be convenient to codify these data
by a single $\q$-divisor, called the {\it branch divisor}
of the orbifold,
$$
\Delta:=\sum (1-\tfrac1{m_j})D_j.
$$

It turns out that the orbifold is
uniquely determined by the pair $(X,\Delta)$.
Slightly inaccurately, we  sometimes  identify the orbifold
with the pair $(X,\Delta)$.

In the cases that we consider $X$ is algebraic,
the $U_i$ are affine, $V_i\cong \c^n$
and the $G_i$ are cyclic, but these special circumstances
are largely unimportant.
\end{defn}

\begin{defn}[Main examples]\label{main.orb.exmp}

Fix (positive) natural numbers $w_1,\dots, w_m$
and consider the $\c^*$-action on $\c^m$ given by
$
\lambda: ( z_1,\dots, z_m)\mapsto (\lambda^{w_1}z_1,\dots,\lambda^{w_m}z_m).
$
Set $W=\gcd(w_1,\dots, w_m)$.
The $W$th roots of unity act trivially  on  $\c^m$,
hence without loss of generality we can replace the action by
$$
\lambda: ( z_1,\dots, z_m)\mapsto 
(\lambda^{w_1/W}z_1,\dots,\lambda^{w_m/W}z_m).
$$
That is, we can and will assume that 
the $w_i$ are relatively prime, i.e. $W=1.$ It is convenient to
 write the $m$-tuple 
$(w_1,\cdots,w_m)$ in vector notation as $\bfw=(w_1,\cdots,w_m)$, 
and to denote the 
$\bbc^*$ action by $\bbc^*(\bfw)$ when we want to specify the action.

We construct an orbifold by considering the quotient 
of $\bbc^m\setminus\{0\}$ by this $\bbc^*$  
action. We write this quotient as 
$\bbp(\bfw)=(\bbc^m\setminus\{0\})/\bbc^*(\bfw).$ 
The orbifold structure is defined as follows.
Set $V_i:=\{(z_1,\cdots,z_m) ~|~z_i=1\}$. Let $G_i\subset \c^*$ be 
the subgroup of $w_i$-th
roots of unity. Note that $V_i$ is invariant under the action of
$G_i$. Set $U_i:=V_i/G_i$.
Note that the $\c^*$-orbits on $(\bbc^m\setminus\{0\})\setminus (z_i=0)$
are in one--to--one correspondence with the points
of $U_i$, thus we indeed have defined charts of an orbifold.
As an algebraic variety this gives the {\it weighted projective space}
 $\bbp(\bfw)$ defined as 
the projective scheme of the graded polynomial ring
 $S(\bfw)=\bbc[z_1,\cdots,z_m]$, where 
$z_i$ has grading or weight $w_i.$
The weight $d$ piece of  $S(\bfw)$, also denoted by
$H^0(\p({\bf w}),d)$, is the vector space of
weighted homogeneous polynomials of {\it weighted degree} $d$.
That is, those that satisfy
$$
f(\lambda^{w_1}z_1,\dots,\lambda^{w_m}z_m)=\lambda^{d}f(z_1,\dots,z_m).
$$
The weighted degree of $f$ is denoted by $w(f)$.

Let $0\in Y\subset \c^m$ be a subvariety with an isolated singularity
at the origin which  is invariant under the given $\c^*$-action.
Similarly, we can construct an orbifold on the quotient 
$(Y\setminus \{0\})/\c^*(\bfw).$ 
As a point set, it is the set of orbits of $\c^*(\bfw)$ on 
$Y\setminus \{0\}$. It's orbifold 
structure is that induced from the orbifold structure on $\bbp(\bfw)$ 
obtained by intersecting 
the orbifold charts described above with $Y.$
In order to simplify notation, we denote it by $Y/\c^*(\bfw)$
or by  $Y/\c^*$ if the weights are clear.

\end{defn}

\begin{defn}
Many definitions concerning orbifolds simplify if we introduce
an open set $U_{ns}\subset X$ which is the complement of the
singular set of $X$ and of the branch divisor.
Thus $U_{ns}$ is smooth and we take $V_{ns}=U_{ns}$.

For the main examples described above
$U_{ns}$  is exactly the set of those orbits
where the stabilizers are trivial.
Every orbit contained in $\c^m\setminus(\prod z_i=0)$ is such.
More generally, a point $(y_1,\dots,y_m)$ corresponds to such an orbit 
if and only if $\gcd\{w_i:y_i\neq 0\}=1$.
\end{defn}

\begin{defn}[Tensors on orbifolds]
\label{diff.orb.defn}

A {\it tensor} $\eta$ on the orbifold $(X,\Delta)$
is a tensor $\eta_{ns}$ on $U_{ns}$
such that for every chart $\phi_i:V_i\to U_i$
the pull back $\phi_i^*\eta_{ns}$ extends to a
tensor  on $V_i$.
In the {\it classical case} the complement of $U_{ns}$
has codimension at least 2, so by Hartogs' theorem
holomorphic tensors on $U_{ns}$
can be identified with holomorphic tensors on the orbifold.
This is not so if there is a branch divisor $\Delta$.
We are especially interested in understanding
the top dimensional holomorphic forms and their tensor powers.

The {\it canonical line bundle} of the orbifold
$K_{X^{orb}}$ is a family of line bundles, one on each chart $V_i$, which is 
the highest exterior power of the holomorphic cotangent bundle
$\Omega^1_{V_i}=T_{V_i}^*$. We would like to study
global sections of powers of $K_{X^{orb}}$.
Let $U_i^{ns}$ denote the smooth part of $U_i$
and $V_i^{ns}:=\phi_i^{-1}U_i^{ns}$.
As shown by (\ref{orbif.defn}.1), 
$K_{V_i}$ is not the pull back of $K_{U_i}$, rather
$$
K_{V_i^{ns}}\cong \phi_i^*K_{U_i^{ns}}(\sum (m_{ij}-1)R_{ij}).
$$
Since $R_{ij}=m_j\phi_i^*D_{ij}$, we obtain, at least formally,  that
$K_{X^{orb}}$ is the pull back of $K_X+\Delta$, rather than the
pull back of $K_X$. The latter of course makes sense only if we define
fractional tensor powers of 
line bundles. Instead of doing it, we state a consequence
of the formula:
\end{defn}

\begin{claim} For $s>0$, global sections of
$K_{X^{orb}}^{\otimes s}$
are those sections of $K_{U_{ns}}^{\otimes s}$ which have an at most
$s(m_i-1)/m_i$-fold pole along the branch divisor $D_i$ for every $i$.
For $s<0$, global sections of
$K_{X^{orb}}^{\otimes s}$
are those sections of $K_{U_{ns}}^{\otimes s}$ which have an at least
$s(m_i-1)/m_i$-fold zero along the branch divisor $D_i$ for every $i$.
\end{claim}

\begin{defn}[Metrics on orbifolds]\label{orb.met.defn}
A {\it Hermitian metric} $h$ on the orbifold $(X,\Delta)$
is  a Hermitian metric $h_{ns}$ on $U_{ns}$
such that for every chart $\phi_i:V_i\to U_i$
the pull back $\phi^*h_{ns}$ extends to a
Hermitian metric on $V_i$.
One can now talk about curvature, K\"ahler metrics, \Ke metrics
on orbifolds. 
\end{defn}

\begin{say}[The  hypersurface case]\label{hypsf.K.say}

We are especially interested in the  case
when $Y\subset \c^m$ is a hypersurface.
It is then the zero set of a polynomial
$F(z_1,\dots,z_m)$ which is equivariant
with respect to the $\c^*$-action. 
$F$ is irreducible since it has an isolated singularity at the origin,
and we  always assume that $F$ is not one of the $z_i$.
Thus $Y\setminus(\prod z_i=0)$ is dense in $Y$.

A differential form on $U_{ns}$ is the same as
a $\c^*$-invariant differential form on $Y_{ns}$
and such a form corresponds to a global
 differential form on $X^{orb}$ iff 
the corresponding $\c^*$-invariant differential form
extends to $Y\setminus\{0\}$.

 The $(m-1)$-forms
$$
\eta_i:=\frac1{\partial F/\partial z_i}
dz_1\wedge \dots\wedge  \widehat{dz_i}\wedge \dots\wedge   dz_m|_Y
$$
satisfy $\eta_i=(-1)^{i-j}\eta_j$ and they glue together to a
global generator $\eta$ of the 
canonical line bundle $K_{Y\setminus\{0\}}$ of $Y\setminus\{0\}$. 
\end{say}

\begin{prop}\label{hypsf.K.prop}
 Assume that $m\geq 3$ and  $s(w(F) -\sum w_i)>0$. Then the 
following three spaces are naturally isomorphic:
\begin{enumerate}
\item Global sections of
$K_{X^{orb}}^{\otimes s}$.
\item 
$\c^*$-invariant global sections of
$K_Y^{\otimes s}$.
\item The space of
weighted homogeneous polynomials of weight $s(w(F) -\sum w_i)$,
modulo multiples of $F$.
\end{enumerate}
\end{prop}

Proof. We have already established that global sections of
$K_{X^{orb}}^{\otimes s}$ can be identified with
$\c^*$-invariant global sections of
$K_{Y\setminus\{0\}}^{\otimes s}$.
If $m\geq 3$ then $Y$ is a hypersurface of dimension 
 $\geq 2$ with an isolated singularity at
the origin, thus normal. Hence 
 global sections of
$K_Y^{\otimes s}$ agree with global sections of
$K_{Y\setminus\{0\}}^{\otimes s}$. This shows the equivalence of (1) and (2).

The $\c^*$-action on $\eta$ has weight $\sum w_i-w(F)$,
thus  $K_Y^{\otimes s}$
is the trivial bundle on $Y$, where the $\c^*$-action has
weight $s(\sum w_i-w(F))$. Its invariant global sections are thus
given by  homogeneous polynomials of weight $s(w(F) -\sum w_i)$
times the generator $\eta$.\qed

In particular, we see that:

\begin{cor}\label{hypsf.K.amp.cor}
Notation as in \S\ref{hypsf.K.say}. 
$K_{X^{orb}}^{-1}$ is ample
iff $w(F)<\sum w_i$.
\end{cor}

\begin{say}[Automorphisms and Deformations]\label{auts.dfes}

If $m\geq 4$ and $Y\subset \c^m$ is a hypersurface,
then by  the Grothendieck--Lefschetz theorem,
every orbifold line bundle on $Y/\c^*$ is the restriction of
an orbifold line bundle on $\c^m/\c^*$ \cite{sga2}.
 This implies that every isomorphism between two
 orbifolds $Y/\c^*(\bfw)$ and $Y'/\c^*(\bfw')$ 
is induced by an
automorphisms of $\c^m$ which commutes with the $\c^*$-actions.
Therefore the weight sequences $\bfw$ and $\bfw'$ are the same 
(up to permutation) and 
every such automorphism $\tau$ has the form
$$
\tau(z_i)=g_i(z_1,\dots,z_m) \qtq{where} w(g_i)=w_i.
\eqno{\ref{auts.dfes}.1}
$$
They form a group $\aut(\c^m,{\bf w})$.
For small values of $t$, maps of the form
$
\tau(z_i)=z_i+tg_i(z_1,\dots,z_m) \qtq{where} w(g_i)=w_i
$
are  automorphisms, hence the dimension of $\aut(\c^m,{\bf w})$
 is $\sum_i \dim H^0(\p({\bf w}),w_i).$
Thus we see that, up to isomorphisms, the orbifolds 
$Y(F)/\c^*$ where $w(F)=d$  form a family
of complex dimension at least
$$
\dim H^0(\p({\bf w}),d)-\sum_i \dim H^0(\p({\bf w}),w_i),
\eqno{\ref{auts.dfes}.2}
$$
and equality holds if the general orbifold in the family has
only finitely many automorphisms.

\end{say}

\begin{say}[Contact structures]\label{contact.say}
 A {\it holomorphic contact structure}
on a complex manifold $M$ of dimension $2n+1$
is a line subbundle $L\subset \Omega^1_M$
such that if $\theta$ is a local section of $L$ then
$\theta\wedge (d\theta)^n$ is nowhere zero.
This forces  an isomorphism $L^{n+1}\cong K_M$.
We would like to derive necessary conditions for
$X^{orb}=Y/\c^*$ to have an orbifold contact structure.

First of all, its dimension has to be odd, so $m=2n+3$
and $n+1$ must divide the canonical class
$K_{X^{orb}}\cong \o(w(F)-\sum w_i)$.
If these conditions are satisfied, then a contact structure
gives a global section of 
$$
\Omega^1_{X^{orb}}\otimes \o\left(\tfrac2{m-1}(-w(F)+\sum w_i)\right).
$$
By pull back, this corresponds  to a
global section of $\Omega^1_{Y\setminus \{0\}}$ on which 
$\c^*$ acts with weight $\tfrac2{m-1}(-w(F)+\sum w_i)$.

Next we claim that every global section of $\Omega^1_{Y\setminus \{0\}}$
lifts to a global section of $\Omega^1_{\c^m}$.
As a preparatory step, it is easy to compute   that
$H^i(\c^m\setminus\{0\}, \o_{\c^m\setminus\{0\}})=0$ for $0<i<m-1$.
(This is precisely the computation done in \cite[III.5.1]{harts}.)
Using the exact sequence
$$
0\to \o_{\c^m\setminus\{0\}}\stackrel{F}{\longrightarrow}
\o_{\c^m\setminus\{0\}}\to \o_{Y\setminus \{0\}}\to 0,
$$
these imply that 
$H^i(Y\setminus\{0\}, \o_{Y\setminus\{0\}})=0$ for $0<i<m-2$.
Next apply the $i=1$ case
to the   co-normal sequence (cf.\ \cite[II.8.12]{harts})
$$
0\to \o_{Y\setminus\{0\}}\stackrel{dF}{\longrightarrow}
\Omega^1_{\c^m\setminus\{0\}}|_{Y\setminus \{0\}} \to
\Omega^1_{Y\setminus \{0\}} \to 0
$$
to conclude that for  $m\geq 4$, every global section 
of $\Omega^1_{Y\setminus \{0\}}$ lifts  to
a global section of $\Omega^1_{\c^m\setminus\{0\}}|_{Y\setminus \{0\}}$.
The latter is the restriction of the free sheaf 
$\Omega^1_{\c^m}|_{Y}$ to $Y\setminus \{0\}$;
hence, we can extend the global sections to
$\Omega^1_{\c^m}|_{Y}$ since $Y$ is normal. Finally these lift to
global sections of $\Omega^1_{\c^m}$ since $\c^m$ is affine.
$\Omega^1_{\c^m}=\sum_i dz_i\o_{\c^m}$, hence there every
$\c^*$-eigenvector has weight at least $\min_i\{w_i\}$.
So we obtain:

\begin{lem}\label{contact.claim}
 The hypersurface $Y/\c^*$ has no 
holomorphic orbifold contact structure if $m\geq 4$ and 
$
\tfrac2{m-1}(-w(F)+\sum w_i)<\min_i\{w_i\}.
$
\end{lem}

This condition is satisfied for all the orbifolds
considered in Theorem~\ref{BP.KE.thm}.
\end{say}

\section{Sasakian-Einstein structures on links}

\begin{say}[Brief review of Sasakian geometry]\label{revSasa} For more 
details see 
\cite{BG00} and references therein.
Roughly speaking a {\it Sasakian structure} on a manifold $M$ is a contact
 metric structure 
$(\xi,\eta,\Phi,g)$ such that the Reeb vector field $\xi$ is a Killing 
vector field of unit length, 
and whose structure transverse to the flow of $\xi$ is K\"ahler. 
Here $\eta$ is a contact 
1-form, $\Phi$ is a $(1,1)$ tensor field which defines a complex 
structure on the contact 
subbundle $\ker~\eta$ which annihilates $\xi,$ and the metric is
 $g=d\eta\circ (\Phi\otimes {\rm 
id})+\eta\otimes \eta.$

We are interested in the case when both $M$ and the leaves of the 
foliation generated by 
$\xi$ are compact. In this case the Sasakian structure is 
called {\it quasi-regular}, and the 
space of leaves $X^{orb}$ is a compact 
K\"ahler orbifold \cite{BG00}. 
 $M$ is the total space 
of a circle {\it orbi-bundle} (also called V-bundle) over $X^{orb}.$
 Moreover, the 2-form $d\eta$ pushes 
down to a K\"ahler form 
$\gro$ on $X^{orb}.$  Now $\gro$ defines an integral class $[\gro]$ of the 
orbifold cohomology 
group $H^2(X^{orb},\bbz)$ which generally is only a rational class in 
the ordinary cohomology 
$H^2(X,\bbq).$

This construction can be inverted in the sense that given a K\"ahler 
form $\gro$ on a 
compact complex orbifold $X^{orb}$ which defines an element 
$[\gro]\in H^2(X^{orb},\bbz)$ one can 
construct a circle  orbi-bundle
on $X^{orb}$ whose orbifold first Chern 
class is $[\gro].$ Then the total 
space $M$ of this orbi-bundle has a natural Sasakian structure 
$(\xi,\eta,\Phi,g)$, where $\eta$ 
is a connection 1-form whose curvature is $\gro.$ The tensor 
field $\Phi$ is obtained by lifting 
the almost complex structure $I$ on $X^{orb}$ to the horizontal 
distribution $\ker~\eta$ and 
requiring that $\Phi$ annihilates $\xi.$ Furthermore, 
the map $(M,g)\ra{1.3} (X^{orb},h)$ is an 
orbifold Riemannian submersion.

The Sasakian structure constructed by the inversion process is not unique. 
One can perform a 
gauge transformation on the connection 1-form $\eta$ and obtain 
a distinct Sasakian 
structure. However, a straightforward curvature computation shows 
that there is a unique 
Sasakian-Einstein metric $g$ with scalar curvature necessarily $2n(2n-1)$
if and only if the K\"ahler metric $h$ is K\"ahler-Einstein with scalar 
curvature $4(n-1)n$, see \cite{Be87,BG00}. 
Hence, the correspondence between orbifold K\"ahler-Einstein metrics
 on $X^{orb}$ with scalar 
curvature $4(n-1)n$ and \Se metrics on $M$ is one-to-one.

\end{say}

\begin{say}[Sasakian structures on links of isolated hypersurface 
singularities]\label{linkdef}

Let $F$ be a weighted homogeneous polynomial as  in 
Definition~\ref{main.orb.exmp}, 
and consider the subvariety $Y:=(F=0)\subset \bbc^{n+1}.$ Suppose further 
that $Y$ has only an 
isolated singularity at the origin. Then the  link 
$L_F=F^{-1}(0)\cap S^{2m-1}$ of $F$ is a 
smooth compact $(m-3)$-connected manifold of dimension $2m-3$ \cite{mil68}.  
So if $m\geq 4$ the manifold $L_F$ is simply connected.
$L_F$ inherits a circle action from the circle subgroup of the $\bbc^*$ group  
described in 
Definition \ref{main.orb.exmp}. 
We denote this circle group by $S_\bfw^1$ to emphasize its 
dependence on the weights.

As noted in \S \ref{revSasa} the K\"ahler structure on $Y/\c^*$
 induces a Sasakian structure 
on the link $L_F$  
such that the infinitesimal generator of the 
weighted circle action defined on $\bbc^{m}$ 
restricts to the Reeb vector field of the Sasakian structure, 
which we denote by $\xi_\bfw.$  
This Sasakian structure $(\xi_\bfw,\eta_\bfw,\Phi_\bfw,g_\bfw)$,  which is 
induced from the {\it weighted Sasakian structure}
 on $S^{2m-1},$ was first noticed by 
Takahashi \cite{Tak78} for Brieskorn manifolds, and is discussed in detail 
in \cite{BG01}.

The quotient space of the link $L_F$ by this circle action is just the
 orbifold 
$X^{orb}=Y/\bbc^*$ 
introduced in Definition \ref{main.orb.exmp}. 
It has a natural K\"ahler structure. In fact, all of this 
fits 
nicely into a commutative diagram \cite{BG01}:
\begin{equation}\label{commdiag}
\begin{matrix}
L_{F} &\ra{2.5}& S^{2m-1}_\bfw&\\
  \decdnar{\pi}&&\decdnar{\phantom{\pi}} &\\
  \hphantom{orb}X^{orb} &\ra{2.5} &\bbp(\bfw),&
\end{matrix}
\end{equation}
where $S^{2m-1}_\bfw$ emphasizes the weighted Sasakian structure 
described for example in 
\cite{BG01}, the horizontal arrows are Sasakian and K\"ahlerian embeddings,
respectively, and the vertical arrows are orbifold Riemannian submersions. 
In particular, the 
Sasakian metric $g$ satisfies $g=\pi^*h+\eta\otimes \eta$, where $h$ is 
the K\"ahler metric on 
$X^{orb}.$

\end{say}

\begin{say}[Isometries of Sasakian structures]\label{isom.sasak}

Let $(X_1^{orb},h_1)$ and $(X_2^{orb},h_2)$ be two \Ke
orbifolds and $M_1$ and $M_2$ the corresponding
\Se manifolds. As
 explained in  \S\ref{revSasa}, 
$M_1$ and $M_2$ are isomorphic as Sasakian structures
iff $(X_1^{orb},h_1)$ and $(X_2^{orb},h_2)$ are
biholomorphically isometric. Here we are interested in understanding 
isometries between  $M_1$ and $M_2$.  As we see, with two classes of 
exceptions,
isometries  automatically preserve the Sasakian structure as well.

The exceptional cases are easy to describe:
\begin{enumerate}
\item   $M_1$ and $M_2$ are both the sphere $S^{2n+1}$ with its round metric. 
 By a theorem of Boothby and Wang,  the corresponding circle action is 
fixed point free \cite{BoWa} with weights 
$(1,\dots,1)$. 
This happens only in the uninteresting case
when $Y\subset \c^m$ is a hyperplane.
\item $M_1$ and $M_2$ have  a {\it 3--Sasakian}
structure. This means that there is a 2-sphere's 
worth of Sasakian structures with a transitive 
action of $SU(2)$ (cf. \cite{BG99} for precise definitions). 
This happens only if the $X_i^{orb}$ admit holomorphic contact orbifold
structures, see \cite{BG97}.
\end{enumerate}

\begin{thm} Let $(X_1^{orb},h_1)$ and $(X_2^{orb},h_2)$ be two \Ke
orbifolds and $M_1$ and $M_2$ the corresponding
\Se manifolds. Assume that
we are not in either of the exceptional cases enumerated above.

 Let $\phi:M_1\to M_2$ be an isometry.
Then there is an isometry $\bar{\phi}:X_1^{orb}\to X_2^{orb}$
 which is either  holomorphic or  anti-holomorphic,
such that the following digram commutes:
\begin{equation*}\label{commdiag2}
\begin{matrix}
M_1 &\fract{\phi}{\ra{2.5}}& M_2&\\
\decdnar{\pi_1} && \decdnar{\pi_2} &\\
 \hphantom{o}X_1^{orb} &\fract{\bar{\phi}}{\ra{2.5}}& 
\hphantom{o}X_2^{orb}&
\end{matrix}
\end{equation*}
Moreover, $\bar{\phi}$ determines $\phi$ up to the $S^1$-action
given by  the Reeb vector field.
\end{thm}

Proof.
Let $\cals_i$  denote the Sasakian structure
 on $M_i$. 
Then $\cals_1$ and $\phi^*\cals_2$ are Sasakian structures 
on $M_1$ sharing the same 
Riemannian metric.  
Since neither $g_1$ nor $g_2$ are of 
constant curvature nor part of a 3-Sasakian structure, the 
proof of Proposition~8.4 of \cite{BGN03} implies that 
either $\phi^*\cals_2=\cals_1$ or
$\phi^*\cals_2=\cals_1^c$ the conjugate Sasakian structure,  
$\cals_1^c:=(-\xi_1,-\eta_1,-\Phi_1,g_1).$ Thus, $\phi$ intertwines 
the foliations and gives rise to an 
orbifold map $\bar{\phi}:X_1^{orb}\ra{1.3} X_2^{orb}$ 
as required.

Conversely, any such  biholomorphism or anti-biholomorphism 
$\bar{\phi}$ lifts to an 
orbi-bundle map $\phi:M_1\ra{1.5} M_2$ uniquely up to the $S^1$-action
given by  the Reeb vector field.
\qed

Putting this together with \S\ref{auts.dfes}  we obtain:

\begin{cor} Let $Y_1\subset \c^m$ (resp.\  $Y_2\subset \c^m$)
be weighted homogeneous hypersurfaces
 with  isolated singularities at the origin
with
weights $\bfw_1$ (resp.\ $\bfw_2$). Assume that
\begin{enumerate}
\item $m\geq 4$.
\item $Y_1,Y_2$ have isolated singularities at the origin.
\item $Y_1/\c^*(\bfw_1)$ and $Y_2/\c^*(\bfw_2)$ both have \Ke metrics.
\item Neither $Y_1/\c^*(\bfw_1)$ nor $Y_2/\c^*(\bfw_2)$ has a
 holomorphic contact 
structure.
\end{enumerate}
Let $(L_1,g_1)$  and $(L_2,g_2)$ be the corresponding
Einstein metrics on the links.
Then
\begin{enumerate}\setcounter{enumi}{4}
\item The connected component of the isometry group of $(L_i,g_i)$ is
the circle $S^1$.
\item $(L_1,g_1)$  and $(L_2,g_2)$ are isometric iff
$\bfw_1=\bfw_2$  (up to permutation) and there is an
automorphism $\tau\in \aut(\c^m,\bfw_1)$ as in
(\ref{auts.dfes}.1)  such that $\tau(Y_1)$ is either $Y_2$
or its conjugate $\bar Y_2$.
\end{enumerate}
\end{cor}

\end{say}

\section{K\"ahler-Einstein metrics on orbifolds}

\begin{say}[Continuity method for finding \Ke metrics]
\label{KE.cont.meth}

Let $(X,\Delta)$ be a compact orbifold of dimension $n$
such that $K_{X^{orb}}^{-1}$ is ample.
The continuity method for finding a \Ke metric
on $(X,\Delta)$ was developed by
\cite{aub, siu1, siu2, tian, nad, dem-koll}.

We start with an arbitrary smooth Hermitian metric
$h_0$ on $K_{X^{orb}}^{-1}$
with positive definite curvature form $\theta_0$.
Choose a K\"ahler metric $\omega_0$ such that
$\ric(\omega_0)=\theta_0$.
Since $\theta_0$ and $\omega_0$ represent the same cohomology class,
there is a $C^{\infty}$ function $f$ such that
$$
\omega_0=\theta_0+\tfrac{i}{2\pi}\partial\bar\partial f.
$$
Our aim is to find a family of functions
$\phi_t$ and numbers $C_t$ for  $t\in [0,1]$, 
normalized by the condition
$\int_X\phi_t\omega_0^n=0$, such that they
satisfy the
Monge--Amp\`ere equation
$$
\log\frac{(\omega_0+\tfrac{i}{2\pi}\partial\bar\partial\phi_t)^n}
{\omega_0^n}+t(\phi_t+f)+C_t=0.
$$
We start with $\phi_0=0, C_0=0$ and if we can reach $t=1$,
we get a \Ke metric
$$
\omega_1=\omega_0+\tfrac{i}{2\pi}\partial\bar\partial\phi_1.
$$
It is easy to see that solvability is an open condition
on $t\in [0,1]$, the hard part is closedness.
It turns out that the critical step is a 0th order estimate.
That is, as the values of $t$ for which the Monge--Amp\`ere equation
is solvable approach a critical value $t_0\in [0,1]$,
a subsequence of the $\phi_t$ converges to 
a function $\phi_{t_0}$ which is 
the sum of a $C^{\infty}$ and of 
a plurisubharmonic function. 
 By \cite{tian} we only
need to prove that
$$
\int_Xe^{-\gamma \phi_{t_0}}\omega_0^n<+\infty
\qtq{for some $\gamma>\tfrac{n}{n+1}$.}
\eqno{\ref{KE.cont.meth}.1}
$$

We view $h_0e^{-\phi_{t_0}}$ as a singular metric on
$K_{X^{orb}}^{-1}$. Its curvature current
$\theta_0+\tfrac{i}{2\pi}\partial\bar\partial\phi_{t_0}$
is easily seen to be semi positive.

The method is thus guaranteed to work
if there is no singular metric with semi positive curvature on 
$K_{X^{orb}}^{-1}$ for which the integral in (\ref{KE.cont.meth}.1)
 is divergent.

A theorem of Demailly and Koll\'ar  establishes
how to approximate a plurisubharmonic function  by
sums of logarithms of absolute values of  holomorphic functions 
\cite{dem-koll}.
This allows us to replace an arbitrary
plurisubharmonic function $\phi_{t_0}$ by $\tfrac1{s}\log|\tau_s|$, where
$\tau_s$ is holomorphic. This gives the following criterion:

\end{say}

\begin{thm}\cite{dem-koll} \label{DK.cond.thm}
Let $X^{orb}$ be a compact, $n$-dimensional orbifold such that
 $K_{X^{orb}}^{-1}$ is ample.
The continuity method produces a
\Ke metric on $X^{orb}$ if the following holds:

There is a $\gamma>\tfrac{n}{n+1}$ such that 
 for every
$s\geq 1$ and for every holomorphic section
$\tau_s\in H^0(X^{orb}, K_{X^{orb}}^{-s})$
the following integral is finite:
$$
\int  |\tau_s|^{-\frac{2\gamma}{s}}\omega_0^n\  <\  +\infty.
$$ 
\end{thm}

For the hypersurface case considered in \S\ref{hypsf.K.say}
we can combine this with the description of
sections of $H^0(X^{orb}, K_{X^{orb}}^{-s})$
given in Proposition~\ref{hypsf.K.prop} to
make the condition even more explicit:

\begin{cor}\label{DK.cond.cor2} Let $Y=(F(z_1,\dots,z_m)=0)$ be as in 
\S\ref{hypsf.K.say}.
Assume that $w(F)<\sum w_i$. 
The continuity method produces a
\Ke metric on $Y/\c^*$ if the following holds:

There is a $\gamma>\tfrac{n}{n+1}$ such that 
for every weighted homogeneous polynomial $g$ of
weighted degree $s(\sum w_i-w(F))$, not identically zero on $Y$, the function 
$$
|g|^{-\gamma/s}\qtq{is locally $L^2$ on $Y\setminus\{0\}$.}
$$
\end{cor}

In general it is not easy to decide if a given function
$|g|^{-c}$ 
is locally $L^2$ or not, but we at least have the following
easy criterion. (See, for instance, \cite[3.14, 3.20]{koll-sings}.)

\begin{lem}\label{mult1.L2.lem}
 Let $M$ be a complex manifold and $h$ a holomorphic function on $M$.
 If $c\cdot \mult_ph<1$ for every $p\in M$ then 
$|h|^{-c}$ is   locally $L^2$.
\end{lem}

For $g$ as in Corollary~\ref{DK.cond.cor2} it 
is relatively esay to estimate the multiplicities of its zeros
via intersection
theory, and we obtain the following generalization of
\cite[Prop.11]{jk1}.

\begin{prop}\label{JK.L2.cond.prop}
  Let $Y=(F(z_1,\dots,z_m)=0)$ be as in \S\ref{hypsf.K.say}.
Assume that the intersections of $Y$ with any number of hyperplanes
$(z_i=0)$ are all smooth outside the origin.
Let $g$ be a weighted homogeneous polynomial and pick $\delta_i>0$.
Then 
$$
|g|^{-c}\prod_i |z_i|^{\delta_i-1}\qtq{is locally $L^2$ on $Y\setminus\{0\}$}
$$
if
$c\cdot w(F)\cdot w(g) <\min_{i,j}\{w_iw_j\}.$
\end{prop}

Proof. The case when every $\delta_i=1$ is \cite[Prop.11]{jk1}
combined with Lemma~\ref{mult1.L2.lem}. These also show
that in our case the $L^2$-condition holds away from the
hyperplanes $(z_i=0)$.

We still need to check the $L^2$ condition along the divisors
$H_i:=(z_i=0)\cap Y\setminus \{0\}$. This is accomplished by
reducing the problem to an analogous problem on $H_i$
and using induction.

In algebraic geometry, this method is called
{\it inversion of adjunction}. Conjectured by Shokurov, 
the following version is due to  Koll\'ar \cite[17.6]{koll-inv}. 
It was observed by \cite{maniv} that it can
also be derived from the $L^2$-extension theorem of 
Ohsawa and Takegoshi \cite{oh-ta}.
See \cite{koll-sings} or \cite{koll-mor} for more detailed expositions.

\begin{thm}[Inversion of adjunction]\label{inv.adj.thm}
Let $M$ be a smooth manifold, $H\subset M$ a smooth divisor
with equation $(h=0)$ and $g$ a holomorphic function on $M$.
Let $g_H$ denote the restriction of $g$ to $H$ and assume that it is not
 identically zero. The following are
equivalent:
\begin{enumerate}
\item $|g|^{-c}|h|^{\delta-1}$ is locally $L^2$ near $H$ for every $\delta>0$.
\item $|g_H|^{-c}$ is locally $L^2$ on $H$.
\end{enumerate}
\end{thm}

\qed

\section{Differential Topology of Links}

In this section we briefly describe the differential topology 
of odd dimensional spheres that 
can 
be realized as links of Brieskorn--Pham singularities
and discuss methods for determining 
their 
diffeomorphism type. 

\begin{say}[The group $bP_{2m}$]
\label{bP2n.defn}
The essential work here is that of 
Kervaire and Milnor \cite{ker-mil} who showed that associated with each
sphere $S^n$ with $n\geq 5$ there is an Abelian group $\Theta_n$ consisting of
equivalence classes of homotopy spheres $S^n$ that are equivalent under 
oriented h-cobordism. By Smale's h-cobordism theorem this implies equivalence 
under oriented diffeomorphism. The group operation on $\Theta_n$ is connected 
sum. $\Theta_n$ has a subgroup $bP_{n+1}$ consisting of equivalence classes of 
those homotopy $n$-spheres which bound parallelizable manifolds $V_{n+1}.$ 
Kervaire and Milnor \cite{ker-mil} proved that $bP_{2k+1}=0$  for $k\geq 1.$ 
Moreover, 
for $m\geq 2,$ $bP_{4m}$ is cyclic of order
$$|bP_{4m}|=2^{2m-2}(2^{2m-1}-1)~\hbox{numerator}~\!\!\biggl(\frac{4B_m}{m}\biggr),$$
where $B_m$ is the $m$-th Bernoulli number. Thus, for example
 $|bP_8|=28, |bP_{12}|=992, 
|bP_{16}|=8128.$ In the first two cases these include all exotic
 spheres; whereas, in the last 
case 
$|bP_{16}|$ is precisely half of the homotopy spheres.

For $bP_{4m+2}$ the situation is 
still not entirely understood. It entails computing the Kervaire invariant, 
which is hard. It is known (see the recent review paper \cite{La00} and
 references 
therein) that  $bP_{4m+2}=0$ or $\bbz_2$ and is $\bbz_2$ if $m\neq 2^i-1$ 
for any $i\geq 3.$ Furthermore, $bP_{4m+2}$ vanishes for $m=1,3,7,$ and $15.$

\end{say}

To a sequence $\bfa=(a_1,\dots,a_m)\in \bbz_+^{m}$ Brieskorn associates a 
graph 
$G(\bfa)$ whose $m$ vertices are labeled by $a_1,\dots,a_m.$ Two vertices 
$a_i$ and $a_j$ are connected if and only if $\gcd(a_i,a_j)>1.$ Let 
$G(\bfa)_{ev}$ 
denote the connected component of $G(\bfa)$ determined by the even integers. 
Note that all even vertices belong to $G(\bfa)_{ev},$
 but $G(\bfa)_{ev}$ may contain odd 
vertices as well.

\begin{thm}\cite{briesk}\label{Brieskorngraph}
The link $L(\bfa)$ (with $m\geq 4$) is homeomorphic to the $(2m-3)$- sphere 
if and only if either of the following 
hold.
\begin{enumerate}
\item $G(\bfa)$ contains at least two isolated points, or 
\item $G(\bfa)$ contains a unique 
 odd isolated point and $G(\bfa)_{ev}$ has an odd number 
of vertices with $\gcd(a_i,a_j)=2$ 
for any distinct $a_i,a_j\in G(\bfa)_{ev}$.   
\end{enumerate}
\end{thm}

\begin{say}[Diffeomorphism types of the links $L(\bfa)$]\label{difftypes}
In order to distinguish the diffeomorphism types of the links $L(\bfa)$
 we need to treat the cases 
$m=2k+1$ and $m=2k$ separately. 

By \cite{ker-mil}, the diffeomorphism type of 
a homotopy sphere $\Sigma$ in  $bP_{4k}$   
 is determined  by the signature 
(modulo  $8|bP_{4k}|$)
 of a parallelizable manifold $M$ 
whose 
boundary is $\Sigma.$ By the Milnor Fibration Theorem \cite{mil68}, 
if $\Sigma=L(\bfa)$, we can take 
$M$  to 
be the Milnor fiber $M^{4k}(\bfa)$ which
 for links of isolated singularities coming from 
weighted 
homogeneous polynomials is diffeomorphic to the hypersurface  $\{\bfz\in 
\bbc^{m} ~|~f(z_1,\dots,z_{m})=1\}.$

Brieskorn shows that the signature of $M^{4k}(\bfa)$ can be 
written combinatorially as
\begin{eqnarray*}
\grt(M^{4k}(\bfa))&=& \#\bigl\{\bfx\in \bbz^{2k+1} 
~|~0<x_i<a_i~\hbox{and}~0<\sum_{i=0}^{2k}\frac{x_i}{a_i} 
<1~\mod 2 \bigr\}\\
& -& \#\bigl\{\bfx\in \bbz^{2k+1} 
~|~0<x_i<a_i~\hbox{and}~1<\sum_{i=0}^{2k}\frac{x_i}{a_i} <2~\mod 2 \bigr\}.
\end{eqnarray*}
Using a formula of Eisenstein, Zagier (cf. \cite{hir71}) has 
rewritten this  as: 
\begin{equation*}\label{zagfor}
\grt(M^{4k}(\bfa))=\frac{(-1)^k}{N}\sum_{j=0}^{N-1}\cot\frac{\pi(2j+1)}{2N} 
\cot\frac{\pi(2j+1)}{2a_1}\cdots \cot\frac{\pi(2j+1)}{2a_{2k+1}},
\eqno{\ref{difftypes}.1}
\end{equation*}
where $N$ is any common multiple of the $a_i$'s. 

For the case of $bP_{4k-2}$ the diffeomorphism type is determined by 
the so-called {\it Arf 
invariant} 
$C(M^{4k-2}(\bfa))\in \{0,1\}.$  
Brieskorn then proves the following:

\begin{prop}\label{BrKer} $C(M^{4k-2}(\bfa))=1$ holds 
if and only if condition 2 of 
Theorem \ref{Brieskorngraph} holds and the one isolated point, 
say $a_0,$ satisfies 
$a_0\equiv 
\pm 
3 \mod 8.$ 
\end{prop}

Following conventional terminology we 
say that $L(\bfa)$ is a {\it Kervaire sphere}
if  $C(M^{4k-2}(\bfa))=1$. A  Kervaire sphere is not 
always 
exotic, but it is exotic when $bP_{4k-2}=\bbz_2.$

\end{say}

\section{Brieskorn--Pham singularities}

\begin{notation}\label{BP.defn} Consider a Brieskorn--Pham singularity
$Y(\bfa):=(\sum_{i=1}^mz_i^{a_i}=0)\subset \c^m.$
Set $C=\lcm(a_i:i=1,\dots,m)$.
$Y(\bfa)$ is invariant under the $\c^*$-action
$$
(z_1,\dots,z_m)\mapsto (\lambda^{C/a_1}z_1,\dots,\lambda^{C/a_m}z_m).
\eqno{\ref{BP.defn}.1}
$$
In the notation of Definition~\ref{main.orb.exmp} we have $w_i=C/a_i$ and
$w=w(F)=C$.
Thus $Y(\bfa)/\c^*$ is a Fano orbifold iff
$1<\sum_{i=1}^m\frac1{a_i}.$

More generally, we consider weighted homogeneous perturbations
$$
Y(\bfa,p):=(\sum_{i=1}^mz_i^{a_i}+p(z_1,\dots,z_m)=0)\subset \c^m,
\qtq{where $w(p)=C$.}
$$
The genericity
condition we need, which is always satisfied by $p\equiv 0$ is:
\begin{enumerate}
\item[(GC)] The intersections of $Y(\bfa,p)$ with any number of hyperplanes
$(z_i=0)$ \\
are all smooth outside the origin.
\end{enumerate}

 In order to formulate the statement,
we  further set
$$
C^j=\lcm(a_i:i\neq j),\quad
b_j=\gcd(a_j,C^j)\qtq{and}  d_j=a_j/b_j.
$$
\end{notation}

\begin{thm}\label{BP.KE.thm}
 The orbifold $Y(\bfa,p)/\c^*$ is Fano and 
has a \Ke metric if it satisfies condition (GC) and
$$
1<\sum_{i=1}^m\frac1{a_i}<
1+\frac{m-1}{m-2}\min_{i,j}\Bigl\{\frac1{a_i}, \frac1{b_ib_j}\Bigr\}.
$$
\end{thm}

Note that if the $a_i$ are pairwise relatively prime then all the $b_i$'s are 1
and we get the simpler bounds
$1<\sum_{i=1}^m\frac1{a_i}<
1+\frac{m-1}{m-2}\min_{i}\{\frac1{a_i}\}.$

Proof. 
By Corollary~\ref{DK.cond.cor2} we need to show that for every $s>0$ and 
for every weighted homogeneous polynomial $g$ of
weighted degree $s(\sum w_i-w(F))=sC(\sum a_i^{-1}-1)$, the
function 
$$
|g|^{-\gamma/s}\qtq{is locally $L^2$ on $Y\setminus\{0\}$.}
$$

Our aim is to reduce this to a problem on a perturbation of the simpler
Brieskorn--Pham  singularity
$Y(\bfb).$

\begin{lem}\label{BP.red.lem}
  Let $g$ be a weighted homogeneous polynomial with
respect to the $\c^*$-action (\ref{BP.defn}.1). Then there is a polynomial $G$
such that
$$
g(z_1,\dots,z_m)=\prod z_i^{e_i}\cdot G(z_1^{d_1},\dots,z_m^{d_m}).
$$
\end{lem}

Proof. Note that $C=d_iC^i$. Thus $d_i$ divides
$C/a_j=d_iC^i/a_j$ for $j\neq i$ but  $C/a_i$
is relatively prime to $d_i$.
Write $g=\prod z_i^{e_i}\cdot g^*$ where $g^*$ is not divisible 
by any $z_i$. Thus $g^*$ has a monomial which does not contain $z_i$,
and so its weight is divisible $d_i$. Thus every time $z_i$ appears,
its exponent must be divisible by $d_i$.\qed
\medskip

Applying this to the defining equation of $Y(\bfa,p)$
we obtain that $p(z_1,\dots,z_m)=p^*(z_1^{d_1},\dots,z_m^{d_m})$
for some polynomial $p^*$.
Set
$$
Y(\bfb,p^*):=(\sum_{i=1}^mx_i^{b_i}+p^*(x_1,\dots,x_m)=0)\subset \c^m.
$$
We have a  map $\pi:Y(\bfa,p)\to Y(\bfb,p^*)$ given by $\pi^*x_i=z_i^{d_i}$
and
$$
|g|=\pi^* \prod |x_i|^{e_i/d_i}\cdot |G(x_1,\dots,x_m)|.
$$
The Jacobian of $\pi$ has $(d_i-1)$-fold zero along
$(z_i=0)$.
Thus
$$
|g|^{-\gamma/s}\qtq{is locally $L^2$ on $Y\setminus\{0\}$}
$$
iff
$$
|G|^{-\gamma/s}\cdot
\prod |x_i|^{-\frac{\gamma e_i}{sd_i}+\frac1{d_i}-1}
\qtq{is locally $L^2$ on $Y^*\setminus\{0\}$.}
$$
The latter condition is guaranteed by Proposition~\ref{JK.L2.cond.prop}.
Indeed,
first we need that each $x_i$ has exponent bigger than $-1$.
This is equivalent to
$e_i< \gamma^{-1} s$. We know that 
$e_iC/a_i\leq \wdeg g=sC(\sum a_i^{-1}-1)$ and so it is enough
to know that
$\sum a_i^{-1}-1<\frac{m-1}{m-2} \frac1{a_i}$.
The latter  is one of our assumptions.

Note that 
$w(\sum x_i^{b_i})=B$, where $B:=\lcm(b_1,\dots,b_m)$ and
going from $G(x_1,\dots,x_n)$
to $G(z_1^{d_1},\dots,z_m^{d_m})$ multiplies the weighted degree by
$C/B$. Thus $w(G)\leq \frac{B}{C}w(g)=sB(\sum \frac1{a_i}-1)$.
Therefore
the last  condition 
$$
c\cdot w(F)\cdot w(g)<\min\{w_i,w_j\}
$$
of Proposition~\ref{JK.L2.cond.prop} becomes
$$
\frac{\gamma}{s}\cdot B\cdot sB\left(\sum \frac1{a_i}-1\right)
<\min\biggl\{\frac{B}{b_i},\frac{B}{b_i}\biggl\}.
$$
After dividing by $B^2$, this becomes our other assumption.\qed

\begin{note}  As algebraic varieties, $Y(\bfa)/\c^*$
is the same as $Y(\bfb)/\c^*$.
In particular, when the $a_i$ are pairwise relatively prime
then all the $b_i=1$ hence, as a variety,
$Y(\bfa)/\c^*\cong \c\p^{m-2}.$ 
The orbifold structure is given by the divisor
$$
\sum_{i=1}^{m-1} \Bigl(1-\frac1{a_i}\Bigr)(y_i=0)+\Bigl(1-\frac1{a_m}\Bigr)(\sum y_i=0).
$$
It would be very interesting to write down the corresponding
\Ke metric explictly. This form would then hopefully 
give a \Ke metric without the 
required upper bound in Theorem~\ref{BP.KE.thm}.
\end{note}

For most cases we get
orbifolds with  finite automorphism groups:

\begin{prop}\label{finite.aut}
 Assume that $m\geq 4$ and all but one of the $a_i$
is at least 3. Then the automorphism group of $\{\sum_i x_i^{a_i}=0\}/\c^*$
is finite.
\end{prop}

Proof. It is enough to prove that there are no continuous families of
isomorphisms of the form
$$
\tau_t(x_i)=x_i+\sum_{j\geq 1}t^jg_{ij}(x_1,\dots,x_m).
$$
By assumption
$$
\sum_i \tau_t(x_i)^{a_i}=\sum_i x_i^{a_i}.
$$
Let $j_0$ be the smallest $j$ such that $g_{ij_0}\neq 0$ for some $i$
and look at the $t^{j_0}$ term in the Taylor expansion of the
left hand side:
$$
\sum_i x_i^{a_i-1}g_{ij_0}(x_1,\dots,x_m)=0.
$$
Note that $w(g_{ij})=w(x_i)$ so as long as $a_i\geq 3$ for all but one $i$,
the terms coming from different values of $i$ do not cancel.
Thus every $g_{ij_0}=0$, a contradiction.\qed

\begin{rem} More generally, the automorphism group
of any $(F=0)/\c^*$ is finite as long as $w_i<\frac12 w(F)$ for all 
but one of the $w_i$s and $(F=0)$ is smooth outside the origin.
Indeed, in this case we would get a relation
$\sum (\partial F/\partial x_i)\cdot g_{ij_0}=0$.
By assumption, the $\partial F/\partial x_i$ form a regular sequence,
and so linear  relationships with polynomial coefficients
 between them are generated by the obvious ``Koszul'' relations
$(\partial F/\partial x_i)\cdot (\partial F/\partial x_j)=
(\partial F/\partial x_j)\cdot(\partial F/\partial x_i)$.
We get a contradiction by degree considerations.
\end{rem}

\section{Numerical examples}

We can summarize our existence result 
for \Se metrics as follows:

\begin{thm}\label{main.numeric.thm} 
For a sequence of natural numbers $\bfa =(a_1,\dots,a_m)$, set
$$
L(\bfa):=\bigl\{\sum_{i=1}^mz_i^{a_i}=0\bigr\}\cap S^{2m-1}(1)\subset \c^m.
$$
\begin{enumerate}
\item $L(\bfa)$ has a \Se metric if
$$
1<\sum_{i=1}^m\frac1{a_i}<
1+\frac{m-1}{m-2}\min_{i,j}\Bigl\{\frac1{a_i}, \frac1{b_ib_j}\Bigr\},
$$
where the $b_i\leq a_i$ are defined before Theorem (\ref{BP.KE.thm}).
\item  $L(\bfa)$ is  homeomorphic to  $S^{2m-3}$ iff the conditions of 
 Theorem (\ref{Brieskorngraph}) are satisfied. 
The diffeomorphism type can  be determined 
as  in Paragraph (\ref{difftypes}).
\item  Given two sequences $\bfa$ and $\bfa'$ satisfying the condition
(\ref{main.numeric.thm}.1), the manifolds $L(\bfa)$ and $L(\bfa')$ are 
isometric iff  $\bfa$ is a permutation of  $\bfa'$.
\end{enumerate}
\end{thm}

Our ultimate aim is to obtain a complete enumeration of all  sequences 
that yield a \Se metric on some  homotopy 
sphere. As a consequence of Theorem (\ref{main.numeric.thm}),
a step toward this goal is  finding all sequences
$a_1,\dots,a_m$ satisfying the inequalities
(\ref{main.numeric.thm}.1).
We  accomplish this in low dimensions via a computer program,
see \cite{bgkt}. Here 
we content ourselves with  obtaining 
some examples which  show the
double exponential growth of the number of cases.

\begin{exmp}  Consider sequences of the form
$\bfa=(2,3,7,m)$. By explicit calculation,
the corresponding link $L(\bfa)$ gives a \Se metric
on $S^5$ if $5\leq m\leq 41$ and $m$ is relatively prime to
at least two of $2,3,7$. This is satisfied in $27$ cases.
\end{exmp}

\begin{exmp}\label{12.dim.on 5}  Among the above cases,  the sequence 
$\bfa=(2,3,7,35)$ is especially noteworthy. 
If $C(u,v)$ is any sufficiently general homogenous septic polynomial, then
the link of
$$
x_1^2+x_2^3+C(x_3,x_4^5)
$$
also gives a \Se metric
on $S^5$.  The relevant automorphism group of $\c^4$ is
$$
(x_1,x_2,x_3,x_4)\mapsto
(x_1,x_2,\alpha_3x_3+\beta x_4^5,\alpha_4x_4).
$$
Hence we get a $2(8-3)=10$ real dimensional family of \Se metrics on $S^5$.

Similarly, the sequence $\bfa=(2,3,7,43,43\cdot31)$ gives a standard 7-sphere
with a $2(43-2)=82$-dimensional family of \Se metrics on $S^7$.
\end{exmp}

\begin{exmp}[Euclid's or Sylvester's sequence]\label{egypt.defn}
(See \cite[Sec.4.3]{gkp} or \cite[A000058]{sequences}.)
 
Consider the sequence defined by the recursion relation 
\begin{equation*}\label{extremseq}
c_{k+1}=c_1\cdots c_k+1=
c_k^2-c_k+1
\end{equation*}
beginning with $c_1=2.$
It starts as
$$
2,3,7,43,1807, 3263443, 10650056950807,...
$$
It is easy to see (cf.\ \cite[4.17]{gkp}) that 
\begin{equation*}\label{recipext}
c_k\geq (1.264)^{2^k}-\tfrac12 \qtq{and}
\sum_{i=1}^m\frac1{c_i}=1-\frac1{c_{m+1}-1}= 1-\frac{1}{c_1\cdots c_m}.
\end{equation*}
\end{exmp}

\begin{exmp} Consider sequences of the form
$\bfa=(a_1=c_1,\dots,a_{m-1}=c_{m-1},a_m)$.
The  troublesome part of the inequalities (\ref{main.numeric.thm}.1)
 is the computation of the $b_i$. However
$b_i\leq a_i$ thus it is sufficient to satisfy the following stronger 
restriction:
$$
1<\sum_{i=1}^m\frac1{a_i}<
1+\frac{m-1}{m-2}\min_{i,j}\Bigl\{\frac1{a_ia_j}\Bigr\}=
1+\frac{m-1}{m-2}\cdot \frac1{a_{m-1}a_m}.
$$
By direct computation this is satisfied if $c_m-c_{m-1}<a_m<c_m$.
At least a third of these numbers are relatively prime to
$a_1=2$ and to $a_2=3$, thus we conclude:
\end{exmp}

\begin{prop} Our methods yield at least 
$\frac13(c_{m-1}-1)\geq \frac13(1.264)^{2^{m-1}}-0.5$
inequivalent families of 
\Se metrics on (standard and exotic) $(2m-3)$-spheres.
\end{prop}

If $2m-3\equiv 1\mod 4$ then by Proposition~\ref{BrKer},
 all these metrics are on the standard
sphere. If $2m-3\equiv 3\mod 4$ then all these metrics are
on on both standard and exotic spheres but we cannot say anything
in general about their distribution.

\begin{exmp} We consider the sequences 
$\bfa=(a_1=c_1,\dots,a_{m-1}=c_{m-1},a_m=(c_{m-1}-2)c_{m-1})$.
Any two of them are relatively prime,
except for $\gcd(a_{m-1},a_m)=c_{m-1}$, 
and the inequalities (\ref{main.numeric.thm}.1)
 are satisfied.

The Brieskorn--Pham polynomial has weighted homogeneous perturbations
$$
x_1^{a_1}+\cdots+x_{m-2}^{a_{m-2}}+G(x_{m-1},x_m^{c_{m-1}-2})
$$
where $G$ is any homogeneous polynomial of degree $c_{m-1}$.
Up to coordinate changes, these form a
family of complex dimension  $c_{m-1}-2$.
Thus we conclude:
\end{exmp}

\begin{prop}\label{big.moduli.exmp} Our methods yield an at least 
$2(c_{m-1}-2)\geq 2((1.264)^{2^{m-1}}-2.5)$-dimensional (real) 
family of pairwise inequivalent  
\Se metrics on some (standard or exotic) $(2m-3)$-sphere.
\end{prop}

As before, if $2m-3\equiv 1\mod 4$ then  these metrics are on the standard
sphere.

\begin{exmp} Consider sequences of the form
$\bfa=(a_1=2c_1,\dots,a_{m-2}=2c_{m-2},\break a_{m-1}=2, a_m)$
where $a_m$ is relatively prime to all the other $a_i$s. 
By easy computation, the condition of
Theorem~\ref{BP.KE.thm} 
 is satisfied if $2c_{m-2}<a_m<2c_{m-1}-2$.

The relatively prime condition is harder to pin down, but it certainly holds
if in addition $a_m$ is a prime number. By the prime number theorem,
the number of primes in the interval $[c_{m-1},2c_{m-1}]$ is about
$$
\frac{c_{m-1}}{\log c_{m-1}}\geq 
\frac{(1.264)^{2^{m-2}}}{2^{m-1}\log 1.264}
\geq (1.264)^{2^{m-1}-4(m-1)},
$$
so it is still doubly exponential in $m$.

By Proposition 
\ref{BrKer}, for even $m$,  $L(\bfa)$  the standard sphere if $a\equiv \pm 1 
\mod 8$ and  the Kervaire sphere if $a\equiv \pm 3 \mod 8.$ 
It is easy to check for all values of $m$ that we get at least one solution
of both types.
Thus we conclude:
\end{exmp}

\begin{prop} Our methods yield a doubly exponential number of
inequivalent families of 
\Se metrics on both the standard and the Kervaire $(4m-3)$-spheres.
\end{prop}

\begin{ack}  We thank Y.-T.\ Siu for answering several questions, 
and E. Thomas for 
helping us with the computer programs.
We received many helpful comments from G.\ Gibbons, D.\ Kotschick, G.\ Tian,
 S.-T.\ Yau and W.\  Ziller.
CPB and KG were partially supported by the NSF under grant number
DMS-0203219 and JK
was partially supported by the NSF under grant number
DMS-0200883. The authors would also like to thank Universit\`a di Roma ``La
Sapienza" for partial support where discussions on this work initiated.
\end{ack}

\providecommand{\bysame}{\leavevmode\hbox to3em{\hrulefill}\thinspace}
\providecommand{\MR}{\relax\ifhmode\unskip\space\fi MR }
\providecommand{\MRhref}[2]{%
  \href{http://www.ams.org/mathscinet-getitem?mr=#1}{#2}
}
\providecommand{\href}[2]{#2}


\end{document}